\documentclass[11pt]{amsart}
\usepackage{mathrsfs}
\usepackage{amssymb}

\usepackage{titletoc}
\usepackage{amscd}
\usepackage{amsmath, amssymb}
\usepackage{amsfonts}
\usepackage[colorlinks,linkcolor=blue,citecolor=blue, pdfstartview=FitH]{hyperref}

\numberwithin{equation}{section}

\theoremstyle{plain}
\newtheorem{thm}{Theorem}[section]
\newtheorem{theorem}[thm]{Theorem}
\newtheorem{lemma}[thm]{Lemma}
\newtheorem{corollary}[thm]{Corollary}
\newtheorem{proposition}[thm]{Proposition}

\theoremstyle{definition}

\newtheorem{question}[thm]{Question}

\newtheorem{remark}[thm]{Remark}

\newtheorem{definition}[thm]{Definition}

\newtheorem{example}[thm]{Example}

\newtheorem{defn-thm}[thm]{Definition-Theorem}

\newcommand{\sO}{{\mathcal O}}

\newcommand{\C}{{\mathbb C}}

\renewcommand{\H}{{\mathbb H}}

\newcommand{\N}{{\mathbb N}}

\newcommand{\R}{{\mathbb R}}

\renewcommand{\S}{{\mathbb S}}

\newcommand{\Z}{{\mathbb Z}}

\newcommand{\qtq}[1]{\quad\mbox{#1}\quad}
\newcommand{\bp}{\bar{\partial}}
\newcommand{\Om}{\Omega}

\newcommand{\ds}{\oplus}

\newcommand{\ts}{\otimes}

\newcommand{\btheorem}{\begin{theorem}}
\newcommand{\etheorem}{\end{theorem}}
\newcommand{\bproposition}{\begin{proposition}}
\newcommand{\eproposition}{\end{proposition}}
\newcommand{\bdefinition}{\begin{definition}}
\newcommand{\edefinition}{\end{definition}}
\newcommand{\bcorollary}{\begin{corollary}}
\newcommand{\ecorollary}{\end{corollary}}
\newcommand{\bproof}{\begin{proof}}
\newcommand{\eproof}{\end{proof}}
\newcommand{\bremark}{\begin{remark}}
\newcommand{\eremark}{\end{remark}}
\newcommand{\eexample}{\end{example}}
\newcommand{\bexample}{\begin{example}}
\newcommand{\la}{\langle}
\newcommand{\elemma}{\end{lemma}}
\newcommand{\blemma}{\begin{lemma}}
\newcommand{\ra}{\rangle}
\newcommand{\sq}{\sqrt{-1}}

\newcommand{\p}{\partial}

\renewcommand{\bar}{\overline}

\renewcommand{\phi}{\varphi}

\newcommand{\ee}{\end{eqnarray*}}
\newcommand{\be}{\begin{eqnarray*}}

\newcommand{\beq}{\begin{equation}}
\newcommand{\eeq}{\end{equation}}

\newcommand{\bd}{\begin{enumerate}}
\newcommand{\ed}{\end{enumerate}}

\renewcommand{\>}{\rightarrow}

\usepackage{fancyhdr}
\renewcommand{\leftmark}{{\scriptsize Levi-Civita Ricci-flat metrics on compact complex manifolds}}

\begin{document}
\title{Levi-Civita Ricci-flat metrics on compact complex manifolds}

 \makeatletter
\let\uppercasenonmath\@gobble
\let\MakeUppercase\relax
\let\scshape\relax
\makeatother

\author{Jie He}
\address{Academy of Mathematics and System Science, Chinese Academy,  Beijing, 100190, China.}
\email{moltke@sina.com}

\author{Kefeng Liu}

\address{Department of Mathematics, University of California at Los Angeles, California 90095}

\email{liu@math.ucla.edu}

\author{Xiaokui Yang}
\date{}
\address{Morningside Center of Mathematics, Academy of Mathematics and\\ Systems Science, Chinese Academy of Sciences, Beijing, 100190, China}
\address{HCMS, CEMS, NCNIS, HLM, UCAS, Academy of Mathematics and\\ Systems Science, Chinese Academy of Sciences, Beijing 100190,
China} \email{\href{mailto:xkyang@amss.ac.cn}{{xkyang@amss.ac.cn}}}

\maketitle

\begin{abstract} In this paper, we study the geometry of compact
complex manifolds with Levi-Civita Ricci-flat metrics and
  prove that compact complex surfaces
admitting Levi-Civita Ricci-flat metrics are K\"ahler Calabi-Yau
surfaces or Hopf surfaces.
\end{abstract}

{\small{\setcounter{tocdepth}{1} \tableofcontents}}

\section{Introduction}

 Einstein manifolds and Einstein metrics are  fundamental topics
 in math physics and differential geometry. It is well-known that,
 the background Riemannian metric of a K\"ahler-Einstein metric is
  Einstein. However, when the ambient complex manifolds are not
  K\"ahler, the relationships between  Hermitian metrics and
  their background Riemannian metrics are complicated and somewhat mysterious. Let $(X,h)$ be a Hermitian manifold and $g$ be the background Riemannian metric.
 On the Hermitian holomorphic tangent bundle $(T^{1,0}X,h)$, there are two
 typical \emph{metric compatible connections}:
\bd \item[(A)] the Chern connection $\nabla$, i.e. the unique
connection $\nabla$ compatible with the Hermitian metric and also
the complex structure $\bp$;

\item[(B)] the Levi-Civita connection $\nabla^{\text{LC}}$, i.e. the
restriction of the \emph{complexified Levi-Civita connection} on
$T_\C X$ to the holomorphic tangent bundle $ T^{1,0}X$. \ed

\noindent The Chern connection is the key object in complex geometry
and the Levi-Civita connection $\nabla^{\text{LC}}$ is a
representative of the Riemannian geometry. It is well-known that
when $(X,h)$ is not K\"ahler, $\nabla$ and $\nabla^{\text{LC}}$ are
not the same. The complex geometry of the Chern connection is
extensively investigated in the literatures by using various methods
(e.g. \cite{CRS,F3, FY, FLY, L, LY12, LY14, LY17, ST2, ST3, STW, T,
TW,TW2, Yang16, Yang17A, Yang17, Yau2}).

 In \cite{LY14}, we introduced the first Aeppli-Chern classes for
 holomorphic line bundles. Let $L\to X$ be a holomorphic line bundle over $X$. The
first Aeppli-Chern class is defined as $$
c^{\text{AC}}_1(L)=\left[-\sq\p\bp\log h\right]_{\text{A}} \in
H^{1,1}_{\text{A}}(X)$$ where $h$ is an arbitrary smooth Hermitian
metric on $L$ and the  Aeppli cohomology is
  $$  H^{p,q}_{\text{A}}(X):=\frac{\text{Ker} \p\bp \cap \Om^{p,q}(X)}{\text{Im} \p \cap \Om^{p,q}(X)+ \text{Im}\bp \cap
  \Om^{p,q}(X)}.$$
 For a complex manifold $X$,
$c_1^{\mathrm{AC}}(X)$ is defined to be
$c_1^{\mathrm{AC}}(K^{-1}_X)$ where $K_X^{-1}$ is the anti-canonical
 bundle of $X$.  Note that, for a Hermitian line bundle $(L,h)$,
the classes $c_1(L)$ and $c_1^{\textrm{AC}}(L)$ have the same
$(1,1)$-form representative $\Theta^h=-\sq\p\bp\log h$ in different
cohomological classes.

 It is well-known that on a Hermitian manifold $(X,\omega)$,
the first Chern Ricci curvature
$$\text{Ric}(\omega)=-\sq\p\bp\log\det(\omega)$$ represents the first
Chern class $c_1(X)$. As an analog, we proved in
\cite[Theorem~1.1]{LY14} that the first Levi-Civita Ricci curvature
$\mathfrak{Ric}(\omega)$ represents the first Aeppli-Chern class
$c_1^{\text{AC}}(X)$.  Hence, it is very natural to study
(non-K\"ahler) Calabi-Yau manifolds by using the first Aeppli-Chern
class $c_1^{\text{AC}}(X)$ and the first Levi-Civita Ricci curvature
$\mathfrak{Ric}(\omega)$.\\

The classification of various Ricci-flat manifolds are important
topics in differential geometry. The following result is fundamental
and well-known, and we refer to the nice paper \cite{T} of V.
Tosatti for
 discussions on Bott-Chern classes and Chern Ricci-flat
metrics.

\btheorem\label{main0} Let $X$ be a compact complex surface. Suppose
$X$ admits a Chern Ricci-flat Hermitian metric $\omega$, i.e.
$\mathrm{Ric}(\omega)=0$. Then $X$ is  minimal and it is exactly one
of the following \bd
\item  Enriques surfaces;\item
bi-elliptic surfaces;\item  K$3$ surfaces;
\item  $2$-tori;
\item {Kodaira surfaces}.
\ed \etheorem

 \noindent As shown in
\cite{LY14}, the Levi-Civita Ricci-flat condition
$\mathfrak{Ric}(\omega)=0$ is equivalent to \beq
\text{Ric}(\omega)=\frac{1}{2}(\p\p^*\omega+\bp\bp^*\omega),\label{MA}\eeq
where $\text{Ric}(\omega)=-\sq\p\bp\log\left(\omega^n\right)$ is the
 Chern Ricci curvature. The equation (\ref{MA}) is not  of
Monge-Amp\`ere type since there are also non-elliptic terms on the
right hand side. As it is well-known,  it is particularly
challenging
 to solve such equations. By using conformal methods, functional analysis and \emph{explicit constructions}, we obtain the following result analogous to Theorem
 \ref{main0}, which also generalizes the previous results in \cite{LY17}.

\btheorem\label{main} Let $X$ be a compact complex surface. Suppose
$X$ admits a Levi-Civita Ricci-flat Hermitian metric $\omega$, i.e.
$\mathfrak{Ric}(\omega)=0$. Then $X$ is  minimal. Moreover, it lies
in one of the following \bd
\item  Enriques surfaces;\item
bi-elliptic surfaces;\item  K$3$ surfaces;
\item  $2$-tori;
\item {Hopf surfaces}.
\ed \etheorem

\bremark Note that, Enriques surfaces, bi-elliptic surfaces, K$3$
surfaces and  $2$-tori are K\"ahler Calabi-Yau surfaces. It is
obvious that the K\"ahler Calabi-Yau metrics are Levi-Civita
Ricci-flat. However, by using Yau's theorem, there exist many
non-K\"ahler Levi-Civita Ricci-flat metrics on each K\"ahler
Calabi-Yau manifold. \eremark

\bremark It is easy to see that, for a Kodaira surface $X$, it has
$c_1^{\mathrm{BC}}(X)=c_1(X)=c^{\mathrm{AC}}_1(X)=0$. By Theorem
\ref{main0}, it has a Chern Ricci-flat metric. However, we can see
from Theorem \ref{main} that it can \textbf{not} support a
Levi-Civita Ricci-flat metric. \eremark

\bremark It is well-known, it is very difficult to write down
explicitly Ricci-flat metrics. In this paper, we obtain Levi-Civita
Ricci-flat metrics on Hopf surfaces by {explicit constructions}. It
also worths to point out that we only construct Levi-Civita
Ricci-flat metrics on Hopf surfaces of class $1$ (see Theorem
\ref{main128}). We conjecture that all Hopf surfaces can support
Levi-Civita Ricci-flat metrics. On the other hand, every Hopf
surface $X$ is a non-K\"ahler Calabi-Yau manifold, i.e. $c_1(X)=0\in
H^2(X,\R)$. However, $X$ can not support a Chern Ricci-flat
Hermitian metric, i.e. a Hermitian metric $\omega$ with
$\text{Ric}(\omega)=0$. In this view point, the existence of
Levi-Civita Ricci-flat Hermitian metrics on Hopf surfaces is quite
exceptional. \eremark

A compact complex surface $X$ is called a Hopf surface if its
universal covering is analytically isomorphic to
$\C^2\setminus\{0\}$. Its fundamental group $\pi_1(X)$ is a finite
extension of an infinite cyclic group generated by a biholomorphic
contraction which takes the form $ (z,w)\>(az,bw+\lambda z^m) $
where $a,b,\lambda \in \C$, $|a|\geq |b|>1$, $m\in \N^*$ and
$\lambda(a-b^m)=0$. There are two different cases:

\bd \item[(I)] the Hopf surface $H_{a,b}$ of class $1$ if
$\lambda=0$;
\item[(II)] the Hopf surface $H_{a,b,\lambda,m}$ of class $0$ if $\lambda\neq
0$ and $a=b^m$. \ed

\noindent Let $H_{a,b}=\C^{2}\setminus\{0\}/\sim$ where
$(z,w)\sim(az,bw)$ and $|a|\geq |b|> 1$. We set $k_1=\log|a|$ and
$k_2=\log|b|$. Define a real smooth function $
\Phi(z,w)=e^{\frac{k_1+k_2}{2\pi}\theta}$ where $\theta(z,w)$ is a
real smooth function defined by $
|z|^2e^{-\frac{k_1\theta}{\pi}}+|w|^2e^{-\frac{k_2\theta}{\pi}}=1.$
 We construct \textbf{explicitly}
Levi-Civita Ricci-flat metrics on $H_{a,b}$ by perturbations and
conformal changes.

\begin{theorem}\label{main128}
On the Hopf surface $H_{a,b}$ of class $1$,  the Hermitian metric
\beq \omega=\Delta^3\left(\frac{1}{
\Phi}\sqrt{-1}\partial\bar\partial\Phi-\frac{1}{2}\sqrt{-1}\partial\bar\partial\log\Phi\right)
\label{LC} \eeq is Levi-Civita Ricci-flat, i.e.
$\mathfrak{Ric}(\omega)=0$,  where
$$  \Delta=\alpha|z|^2
\Phi^{-\alpha}+(2-\alpha)|w|^2\Phi^{\alpha-2} \qtq{and}
\alpha=\frac{2k_1}{k_1+k_2}.$$
\end{theorem}

\bremark  If $a=b$, $H_{a,a}$ is exactly the usual diagonal Hopf
surface. In this case, the Levi-Civita Ricc-flat metric constructed
in Theorem \ref{main128} is the same as that constructed in
\cite[Theorem~6.2]{LY14} or \cite[Theorem~7.3]{LY17}. \eremark

It is well-known that, on a compact K\"ahler manifold $X$, the
K\"ahler Ricci-flat metrics are all Einstein flat metrics.

\begin{question} On a compact complex manifold $X$, does there exist some Levi-Civita
Ricci-flat (non-K\"ahler) Hermitian metric such that the background
Riemannian metric is Einstein?
\end{question}

\noindent\textbf{Acknowledgement.} The second author would like to
thank  Valentino Tosatti for many useful comments and suggestions.
 This work
was partially supported   by China's Recruitment
 Program of Global Experts and  NSFC 11688101.

\vskip 2\baselineskip

\section{Preliminaries}

\subsection{Chern connection on complex manifolds}
 Let $(X,\omega_g)$ be a compact Hermitian manifold. There exists a unique connection $\nabla$ on  the holomorphic tangent bundle $T^{1,0}X$
  which is compatible with the Hermitian metric and also the complex structure of $X$. This connection $\nabla$ is called the Chern connection.  The Chern
  connection $\nabla$
 on $(T^{1,0}X,\omega_g)$ has  curvature components \beq R_{i\bar
j k\bar \ell}=-\frac{\p^2g_{k\bar \ell}}{\p z^i\p\bar z^j}+g^{p\bar
q}\frac{\p g_{k\bar q}}{\p z^i}\frac{\p g_{p\bar \ell}}{\p\bar
z^j}.\eeq  The (first) Chern-Ricci form $\text{Ric}(\omega_g)$ of
$(X,\omega_g)$ has components
$$R_{i\bar j}=g^{k\bar \ell}R_{i\bar jk\bar \ell}=-\frac{\p^2\log\det(g)}{\p z^i\p\bar z^j}$$
which also represents the first Chern class $c_1(X)$ of the complex
manifold $X$. The \emph{Chern scalar curvature} $s_\mathrm{C}$ of
$(X,\omega_g)$ is given by \beq
s_\mathrm{C}=\text{tr}_{\omega_g}\text{Ric}(\omega_g)=g^{i\bar j}
R_{i\bar j}. \eeq The total \emph{Chern scalar curvature} of
$\omega_g$ is \beq \int_X s_\mathrm{C} \omega_g^n=n
\int\text{Ric}(\omega_g)\wedge \omega_g^{n-1},\eeq where $n$ is the
complex dimension of $X$.

\subsection{Bott-Chern classes and Aeppli classes}
The Bott-Chern
  cohomology and the Aeppli cohomology on a compact complex manifold
  $X$ are given by
  $$ H^{p,q}_{\mathrm{BC}}(X):= \frac{\text{Ker} d \cap \Om^{p,q}(X)}{\text{Im} \p\bp \cap \Om^{p,q}(X)}\qtq{and} H^{p,q}_{\mathrm{A}}(X):=\frac{\text{Ker} \p\bp \cap \Om^{p,q}(X)}{\text{Im} \p \cap \Om^{p,q}(X)+ \text{Im}\bp \cap
  \Om^{p,q}(X)}.$$

\noindent Let $\mathrm{Pic}(X)$ be the set of holomorphic line
bundles over $X$. As similar as the first Chern class map
$c_1:\mathrm{Pic}(X)\>H^{1,1}_{\bp}(X)$, there is a \emph{first
Aeppli-Chern class} map \beq c_1^{\mathrm{AC}}:\mathrm{Pic}(X)\>
H^{1,1}_{\mathrm{A}}(X).\eeq Given any holomorphic line bundle $L\to
X$ and any Hermitian metric $h$ on $L$, its curvature form
$\Theta_h$ is locally given by $-\sq\p\bp\log h$.  We define
$c_1^{\mathrm{AC}}(L)$ to be the class of $\Theta_h$ in
$H^{1,1}_{\mathrm{A}}(X)$.  For a complex manifold $X$,
$c_1^{\mathrm{AC}}(X)$ is defined to be
$c_1^{\mathrm{AC}}(K^{-1}_X)$ where $K_X^{-1}$ is the anti-canonical
line bundle. The first Bott-Chern class $c_1^{\mathrm{BC}}(X)$ can
be defined similary.

\subsection{The Levi-Civita connection on the holomorphic tangent
bundle}

Let's recall some elementary settings (e.g. \cite[Section~2]{LY14}).
Let $(M, g, \nabla)$ be a $2n$-dimensional Riemannian manifold with
the Levi-Civita connection $\nabla$. The tangent bundle of $M$ is
also denoted by $T_\R M$.  Let $T_\C M=T_\R M\ts \C$ be the
complexification. We can extend the metric $g$ and the Levi-Civita
connection $\nabla$ to $T_{\C}M$ in the $\C$-linear way. Let
$(M,g,J)$ be an almost Hermitian manifold, i.e., $J:T_\R M\>T_\R M$
with $J^2=-1$, and
 for any $X,Y\in T_\R M$, $ g(JX,JY)=g(X,Y)$. The Nijenhuis tensor $N_J:\Gamma(M,T_\R M)\times \Gamma(M,T_\R
M)\>\Gamma(M,T_\R M)$ is defined as $$
N_J(X,Y)=[X,Y]+J[JX,Y]+J[X,JY]-[JX,JY].$$ The almost complex
structure $J$ is called \emph{integrable} if $N_J\equiv 0$ and then
we call $(M,g,J)$ a Hermitian manifold. We can also extend $J$ to
$T_\C M$ in the $\C$-linear way. Hence for any $X,Y\in T_\C M$, we
still have $ g(JX,JY)=g(X,Y).$
 By Newlander-Nirenberg's
theorem, there exists a real coordinate system $\{x^i,x^I\}$ such
that $z^i=x^i+\sq x^I$ are  local holomorphic coordinates on $M$.
Moreover, we have $T_\C M=T^{1,0}M\ds T^{0,1}M$ where
$$T^{1,0}M=\text{span}_\C\left\{\frac{\p}{\p z^1},\cdots,\frac{\p}{\p
z^n}\right\}\qtq{and}T^{0,1}M=\text{span}_\C\left\{\frac{\p}{\p \bar
z^1},\cdots,\frac{\p}{\p \bar z^n}\right\}.$$ Since $T^{1,0}M$ is a
subbundle of $T_{\C}M$, there is an induced connection
$\nabla^{\mathrm{LC}}$ on the holomorphic tangent bundle $T^{1,0}M$
given by \beq \nabla^{\text{LC}}=\pi\circ\nabla:
\Gamma(M,T^{1,0}M)\stackrel{\nabla}{\rightarrow}\Gamma(M, T_{\C}M\ts
T_{\C}M)\stackrel{\pi}{\rightarrow}\Gamma(M,T_{\C}M\ts T^{1,0}M).
\eeq Let $h=(h_{i\bar j})$ be the corresponding Hermitian metric on
$T^{1,0}M$ induced by $(M,g, J)$. It is obvious that
$\nabla^{\text{LC}}$ is a metric compatible connection on the
Hermitian holomorphic vector bundle $(T^{1,0}M, h)$, and we call
$\nabla^{\text{LC}}$ the \emph{Levi-Civita connection} on the
complex manifold $M$. It is obvious that, $\nabla^{\text{LC}}$ is
determined by the following relations \beq
\nabla^{\text{LC}}_{\frac{\p}{\p z^i}}\frac{\p}{\p
z^k}:=\Gamma_{ik}^p\frac{\p}{\p z^p} \qtq{and}
\nabla^{\mathrm{LC}}_{\frac{\p}{\p \bar z^j}}\frac{\p}{\p
z^k}:=\Gamma_{\bar jk}^p\frac{\p}{\p z^p} \eeq where \beq
\Gamma_{ij}^k=\frac{1}{2}h^{k\bar \ell}\left(\frac{\p h_{j\bar
\ell}}{\p z^i}+\frac{\p h_{i\bar \ell}}{\p z^j}\right), \qtq{and}
\Gamma_{\bar i j}^k=\frac{1}{2} h^{k\bar \ell}\left( \frac{\p
h_{j\bar\ell}}{\p\bar z^i}-\frac{\p h_{j\bar i}}{\p\bar
z^\ell}\right).\eeq The curvature  tensor $\mathfrak{R}\in
\Gamma(M,\Lambda^2 T_{\C}M\ts T^{*1,0}M\ts T^{1,0}M)$ of
$\nabla^{\mathrm{LC}}$ is given by $$ \mathfrak{R}(X,Y)s
=\nabla^{\mathrm{LC}}_{X}\nabla^{\mathrm{LC}}_Ys-\nabla^{\mathrm{LC}}_Y\nabla^{\mathrm{LC}}_Xs-\nabla^{\mathrm{LC}}_{[X,Y]}s$$
for any $X,Y\in T_{\C}M$ and $s\in T^{1,0}M$. A straightforward
computation shows that the curvature tensor $\mathfrak{R}$ has
$(1,1)$ components \beq \mathfrak{R}_{i\bar
jk}^{\ell}=-\left(\frac{\p \Gamma^{\ell}_{ik}}{\p \bar z^j}-\frac{\p
\Gamma^{\ell}_{\bar jk}}{\p z^i}+\Gamma_{ ik}^{s}\Gamma^{\ell}_{\bar
js}-\Gamma_{ \bar jk}^{s}\Gamma^{\ell}_{s i}\right).\label{curv}\eeq

\noindent The \emph{(first) Levi-Civita Ricci curvature}
$\mathfrak{Ric}(\omega_h)$  of $\left(T^{1,0}M, \omega_h,
\nabla^{\mathrm{LC}}\right)$ is
 \beq \mathfrak{Ric}(\omega_h)=\sq\mathfrak
R^{(1)}_{i\bar j}dz^i\wedge d\bar z^j \qtq{with}
\mathfrak{R}^{(1)}_{i\bar j}=\mathfrak{R}_{i\bar j k}^{k}. \eeq  The
\emph{Levi-Civita scalar curvature} $ s_{\text{LC}}$ of
$\nabla^{\mathrm{LC}}$ on $T^{1,0}M$ is  \beq s_{\text{LC}}=h^{i\bar
j}h^{k\bar \ell}\mathfrak{R}_{i\bar j k\bar \ell}. \eeq

\subsection{Special manifolds} Let $X$ be a compact complex
manifold with complex dimension $n\geq 2$. A Hermitian metric
$\omega_g$ is called a Gauduchon metric if $\p\bp\omega_g^{n-1}=0$.
It is proved by Gauduchon (\cite{Ga2}) that, in the conformal class
of each Hermitian metric, there exists a unique Gauduchon metric (up
to  scaling).  A Hermitian metric $\omega_g$ is called a balanced
metric if $d\omega_g^{n-1}=0$ or equivalently $d^*\omega_g=0$. On a
compact complex surface, a balanced metric is also K\"ahler, i.e.
$d\omega_g=0$. It is well-known many Hermitian manifolds can not
support balanced metrics, e.g. Hopf surface $\S^{3}\times \S^1$. It
is also obvious that balanced metrics are Gauduchon.

\vskip 2\baselineskip
\section{Geometry of the  Levi-Civita connections}

\subsection{Some computational formulas.} In this
subsection, we recall some elementary and well-known computational
lemmas on Hermitian manifolds.

\blemma\label{adjoint} Let $(X,\omega)$ be a compact Hermitian
manifold and $\omega=\sq h_{i\bar j}dz^i\wedge d\bar z^j$.
 \beq\p^*\omega=-\sq
\Lambda\left(\bp\omega\right)=-2\sq \Gamma_{\bar j k}^kd\bar z^j
\qtq{and} \bp^*\omega=\sq \Lambda\left(\p\omega\right)=2\sq
\bar{\Gamma_{\bar ik}^k}dz^i.\label{key1}\eeq \elemma \bproof By the
well-known Bochner formula (e.g. \cite{LY12}), $$[\bp^*,L]=\sq
\left(\p+\tau\right)$$ where $\tau=[\Lambda,\p\omega]$, we see
$\bp^*\omega=\sq \Lambda\left(\p\omega\right)=2\sq \bar{\Gamma_{\bar
ik}^k}dz^i$. \eproof

\noindent Let $T$ be the torsion tensor of the Hermitian metric
$\omega$, i.e. \beq T_{ij}^k=h^{k\bar\ell}\left(\frac{\p h_{j\bar
\ell}}{\p z^i}-\frac{\p h_{i\bar \ell}}{\p z^j}\right).\eeq

 \bcorollary\cite[Corollary~4.2]{LY17}\label{A} Let
$(X,\omega)$ be a compact Hermitian manifold. Let $s$ be the
Riemannian scalar curvature of the background Riemannian metric
$\omega$. Then \beq
s=2s_{\mathrm{C}}+\left(\la\p\p^*\omega+\bp\bp^*\omega,\omega\ra-2|\p^*\omega|^2\right)-\frac{1}{2}|T|^2,\label{key-1}\eeq
\beq
s_{\mathrm{LC}}=s_{\mathrm{C}}-\frac{1}{2}\la\p\p^*\omega+\bp\bp^*\omega,\omega\ra=s_{\mathrm{C}}-\la\p\p^*\omega,\omega\ra.\label{010}\eeq

 \ecorollary

\subsection{The first Aeppli-Chern class and Levi-Civita connections.}
The following result is obtained in  \cite[Theorem~1.2]{LY14}(see
also \cite[Theorem~4.1]{LY17}).

\btheorem\label{main-1} Let $(X,\omega)$ be a compact Hermitian
manifold. Then the first Levi-Civita Ricci form
$\mathfrak{Ric}(\omega)$ represents the first Aeppli-Chern class $
c_1^{\mathrm{AC}}(X)$ in $H^{1,1}_{\mathrm{A}}(X)$. Moreover, we
have the Ricci curvature relation \beq
\mathfrak{Ric}(\omega)=\emph{Ric}(\omega)-\frac{1}{2}(\p\p^*\omega+\bp\bp^*\omega).\label{key}
\eeq
 \etheorem

%

\blemma\label{conformalchange} Let $(X,\omega)$ be a compact
Hermitian manifold with complex dimension $n$. Suppose $f\in
C^{\infty}(X,\R)$ and $\omega_f=e^f\omega$. Then we have \beq
 \bp^*_f\omega_f=\bp^*\omega+\sq(n-1)\p f\qtq{and} \bp\bp^*_f\omega_f= \bp\bp^*\omega-\sq(n-1)\p\bp f.\label{key12}\eeq
where $\bp^*, \bp^*_f$ are the adjoint operators with respect to the
metric $\omega$ and $\omega_f$ respectively. In particular, we have
$$\mathfrak{Ric}(e^f\omega)=\mathfrak{Ric}(\omega)-\sq\p\bp f.$$
\elemma

\bdefinition The Kodaira dimension $\kappa(L)$ of a line bundle $L$
is defined to be
$$\kappa(L):=\limsup_{m\>+\infty} \frac{\log \dim_\C
H^0(X,L^{\ts m})}{\log m}$$ and the \emph{Kodaira dimension}
$\kappa(X)$ of $X$ is defined as $ \kappa(X):=\kappa(K_X)$ where the
logarithm of zero is defined to be $-\infty$.

\edefinition

\btheorem\label{scalarvanishing} Let $(X,\omega)$ be a compact
Hermitian manifold of complex dimension. Let $\omega_f=e^{f}\omega$
be the Gauduchon metric
 in the conformal class of $\omega$. Then we have
\beq \int_X s_{{f}}\cdot
\omega_f^n=n\int_X\mathrm{Ric}(\omega_f)\wedge \omega_f^{n-1}=\int_X
e^{(n-1)f}\cdot s_{\mathrm{LC}}\cdot
\omega^{n}+n\|\bp^*_f\omega_f\|^2_{\omega_f}.\label{Gauduchon}\eeq

  \etheorem

\bproof By Lemma \ref{conformalchange} and Theorem \ref{main-1}, we
have
\be\text{Ric}(\omega_f)-\frac{\p\p^*_f\omega_f+\bp\bp_f^*\omega_f}{2}&=&\text{Ric}(\omega)-\frac{\p\p^*\omega+\bp\bp^*\omega}{2}-\sq\p\bp
f\\
&=&\mathfrak{Ric}(\omega)-\sq\p\bp f. \ee

\noindent Moreover, we have
\begin{eqnarray} \int_X\text{Ric}(\omega_f)\wedge
\omega_f^{n-1}\nonumber&=&\int_X\left(\mathfrak{Ric}(\omega)-\sq\p\bp
f+\frac{\p\p^*_f\omega_f+\bp\bp^*\omega_f}{2}\right)\wedge
\omega_f^{n-1}\\
\nonumber&=&\int_X\mathfrak{Ric}(\omega)\wedge \omega_f^{n-1}+
\frac{1}{2}\left(\|\bp^*_f\omega_f\|^2_{\omega_f}+\|\p^*_f\omega_f\|^2_{\omega_f}\right)\\
&=&\frac{1}{n}\int_X e^{(n-1)f}\cdot s_{\mathrm{LC}}\cdot
\omega^{n}+\|\bp^*_f\omega_f\|^2_{\omega_f}.\end{eqnarray} \eproof

\btheorem\label{Kodaira} Let $X$ be a compact complex manifold.
Suppose $\omega$ is a  Hermitian metric with $s_{\mathrm{LC}}\geq
0$. Then either \bd
\item $\kappa(X)=-\infty$; or\item  $\kappa(X)=0$ and $(X,\omega)$
is conformally balanced with $K_X$ a holomorphic torsion, i.e.
$K_X^{\ts m}=\sO_X$ for some $m\in \Z^+$.\ed \etheorem

\bproof Let $\omega_f=e^{f}\omega$ be the Gauduchon metric in the
conformal class of $\omega$. Then by formula (\ref{Gauduchon}), the
total Chern scalar curvature of $\omega_f$ is \beq \int_X
s_{{f}}\cdot \omega_f^n=n\int_X\text{Ric}(\omega_f)\wedge
\omega_f^{n-1}\geq
n\|\bp_f^*\omega_f\|_{\omega_f}^2\label{LCflat}\eeq since the
Levi-Civita scalar curvature $s_{\mathrm{LC}}\geq 0$. Suppose
$\bp^*_f\omega_f\neq 0$,  then
$$\int_X s_{{f}}\cdot
\omega_f^n>0.$$ By \cite[Corollary~3.3]{Yang17}, we have
$\kappa(X)=-\infty$. On the other hand, if $\bp^*_f\omega_f=0$, i.e.
$(X,\omega)$ is conformally balanced. Then the total Chern scalar
curvature of the Gauduchon metric
$$\int_X s_{{f}}\cdot
\omega_f^n\geq 0.$$ Then by \cite[Theorem~1.4]{Yang17}, we have
$\kappa(X)=-\infty$ or $\kappa(X)=0$, and when $\kappa(X)=0$, $K_X$
is a holomorphic torsion.
 \eproof

\vskip 2\baselineskip

\section{Compact complex manifolds with Levi-Civita Ricci-flat metrics}

 Let's recall that, a
Levi-Civita Ricci-flat metric is a Hermitian metric satisfying
$\mathfrak{Ric}(\omega)=0$, or equivalently, by formula
(\ref{key})\beq
\mathrm{Ric}(\omega)=\frac{\p\p^*\omega+\bp\bp^*\omega}{2}.\eeq It
is easy to see that  \bcorollary\label{c1}  Let $X$ be a compact
complex manifold. Then
$$c_1^{\mathrm{BC}}(X)=0 \Longrightarrow c_1(X)=0 \Longrightarrow c_1^{\mathrm{AC}}(X)=0.$$
 \ecorollary

\noindent The first obstruction for the existence of Levi-Civita
Ricci-flat Hermitian metric is the top first Chern number:

\bcorollary\label{top} Suppose $c_1^{\mathrm{AC}}(X)=0$, then the
top intersection number $c_1^{n}(X)=0$. In particular, if $X$ has a
Levi-Civita Ricci-flat Hermitian metric $\omega$, then
$c_1^{n}(X)=0$. \ecorollary

\bproof By definition,  if $c_1^{\mathrm{AC}}(X)=0$, then
$$\text{Ric}(\omega)=\bp A+\p B$$
where $A$ is a $(1,0)$-form and $B$ is a $(0,1)$-form. Hence
$$c_1^n(X)=\int_X (\text{Ric}(\omega))^n=\int_X (\text{Ric}(\omega))^{n-1}\wedge (\bp A+\p B)=0$$
since $\text{Ric}(\omega)$ is both $\p$ and $\bp$-closed.
 \eproof

\btheorem\label{no1} Let $X$ be a compact complex manifold with
$\kappa(X)=-\infty$.  If $X$ has a  Levi-Civita Ricci-flat Hermitian
metric $\omega$, then $X$ must be a non-K\"ahler manifold. \etheorem

\bproof Let $\omega$ be a Hermitian metric with
$\mathfrak{Ric}(\omega)=0$. By formula (\ref{key}), we have \beq
\text{Ric}(\omega)=\frac{\p\p^*\omega+\bp\bp^*\omega}{2}. \eeq Note
that $\text{Ric}(\omega)$ is $\p$-closed and $\bp$-closed, and so we
have
$$\bp\p\p^*\omega=0$$
 Suppose $X$ is a K\"ahler manifold, then  by $\p\bp$-Lemma on $X$,  the $\bp$-closed and $\p$-exact $(1,1)$-form
$\p\p^*\omega$ is $\p\bp$-exact, i.e. there exists a smooth function
$\phi$ such that
$$\p\p^*\omega=\p\bp\phi.$$
Therefore, $$\text{Ric}(\omega)=\sq\p\bp F$$ where
$F=-\frac{\bar\phi-\phi}{2\sq}\in C^\infty(X,\R)$. It is obvious
that the Hermitian metric $e^{\frac{F}{n}}\omega$ is Chern
Ricci-flat, i.e. $c_1^{\mathrm{BC}}(X)=c_1(X)=0$ is unitary flat.
Hence $X$ is a K\"ahler Calabi-Yau manifold. In particular,
$\kappa(X)=0$. This is a contradiction.
 \eproof

\bremark On a compact K\"ahler Calabi-Yau manifold $X$ with $\dim_\C
X=n\geq 2$, the Levi-Civita Ricci-flat metrics are not necessarily
K\"ahler. Indeed, let $\omega_{\mathrm{CY}}$ be a Calabi-Yau
K\"ahler metric on $X$.
 Then for any non constant smooth function $f\in C^\infty(X,\R)$, we can construct a non-K\"ahler
 Levi-Civita Ricci-flat metric. By Yau's  theorem, there
 exists a K\"ahler metric $\omega_0$ such that
 $$\omega_0^n=e^{-f}\omega^n_{\mathrm{CY}}.$$
Let $\omega=e^f\omega_0$. Then $\omega$ is a non-K\"ahler metric
with Levi-Civita Ricci-flat curvature. Indeed, \be
\mathfrak{Ric}(\omega)&=&\text{Ric}(\omega)-\frac{\p\p^*\omega+\bp\bp^*\omega}{2}\\
&=&\text{Ric}(\omega_0)-n\sq\p\bp
f-\frac{\p\bp_0^*\omega_0+\bp\bp^*_0\omega_0}{2}+(n-1)\sq\p\bp f\\
&=&\text{Ric}(\omega_{\mathrm{CY}})+\sq\p\bp f-n\sq\p\bp
f-\frac{\p\bp_0^*\omega_0+\bp\bp^*_0\omega_0}{2}+(n-1)\sq\p\bp f \\
&=&0,\ee where we use Lemma \ref{conformalchange} in the second
identity.

\eremark

Theorem \ref{no1} has the following variant: \bcorollary \label{no}
Let $X$ be a compact K\"ahler manifold. If $\kappa(X)=-\infty$, then
$X$ has no Levi-Civita Ricci-flat Hermitian metric. \ecorollary

 As an application, we obtain

 \btheorem\label{greater} Let $X$ be a compact complex surface with
$\kappa(X)\geq 0$. Suppose $X$ admits a Hermitian metric with
$s_{\mathrm{LC}}\geq 0$. Then $X$ is a minimal K\"ahler surface of
Calabi-Yau type, i.e. $X$ is exactly one of the following \bd\item
an Enriques surface;\item a bi-elliptic surface;\item a K$3$
surface;
\item a torus.\ed \etheorem

\bproof By Theorem \ref{Kodaira}, we know $\kappa(X)=0$ and
 the canonical line bundle $K_X$ is a
holomorphic torsion, i.e. $K_X^{\ts m}=\sO_X$ for some $m\in \Z^+$.
Since $\dim X=2$, by Theorem \ref{Kodaira} again,  $X$ is a balanced
surface and so it is K\"ahler. It is easy to see that, $X$ is
minimal. According to the Kodaira-Enriques' classification, $X$ is
either an Enriques surface, a bi-elliptic surface, a K$3$ surface or
a torus. All these surfaces are K\"ahler surfaces of Calabi-Yau
type, and all K\"ahler Calabi-Yau metrics are Levi-Civita
Ricci-flat. \eproof

\vskip 2\baselineskip

\section{The proof of Theorem \ref{main}}

In this section, we prove Theorem \ref{main}, i.e.

\btheorem\label{main22} Let $X$ be a compact complex surface.
Suppose $X$ admits a Levi-Civita Ricci-flat Hermitian metric
$\omega$. Then $X$ is  minimal.  Moreover, it lies in one of the
following \bd
\item  Enriques surfaces;\item
bi-elliptic surfaces;\item  K$3$ surfaces;
\item  $2$-tori;
\item  Hopf surfaces.
\ed \etheorem

\bproof Let $(X,\omega)$ be a compact complex surface with
Levi-Civita Ricci-flat metric $\omega$. Then we have
$c_1^{\mathrm{AC}}(X)=0$ and $s_{\mathrm{LC}}=0$. By Theorem
\ref{Kodaira}, $\kappa(X)=-\infty$ or $\kappa(X)=0$. We shall show
$X$ is a minimal surface.

Suppose $\kappa(X)=0$, by Theorem \ref{greater}, we know $X$ is a
 minimal K\"ahler Calabi-Yau  surface,  i.e. $X$ is exactly one of the following \bd\item a
Enriques surface;\item a bi-elliptic surface;\item a K$3$ surface;
\item a torus.\ed

 Suppose  $\kappa(X)=-\infty$. Let $X_{\mathrm{min}}$ be the minimal
model of $X$, then $X_{\mathrm{min}}$ lies in one of the following
classes:

 \bd \item minimal rational surfaces;

 \item ruled surfaces of genus $g\geq 1$;

 \item  surface of class $\mathrm{VII}_0$.
 \ed
If $X_{\mathrm{min}}$ is in $(1)$ or $(2)$, we know $X$ is
projective. Since $\kappa(X)=-\infty$,  by Corollary \ref{no}, $X$
has no Levi-Civita Ricci-flat metric. Hence $X_{\mathrm{min}}$ is
not in $(1)$ or $(2)$.  Suppose $X_{\mathrm{min}}$ lies in $(3)$,
i.e.  of class $\mathrm{VII}_0$. A class $\mathrm{VII}_0$ surface is
a minimal compact complex surface with $b_1=1$ and
$\kappa(X)=-\infty$. It is well-known that the first Betti number
$b_1$ of compact complex surfaces  are invariant under blowing-ups,
i.e. $b_1(X)=1$. By \cite[Theorem ~2.7 on p.139]{BHPV}, we know
$$b_1(X)=h^{1,0}(X)+h^{0,1}(X), \qtq{and} h^{1,0}(X)\leq
h^{0,1}(X)$$ hence $h^{0,1}(X)=1$. Since
$\kappa(X)=\kappa(X_{\mathrm{min}})=-\infty$, by Serre duality, we
have
$$h^{0,2}(X)=h^{2,0}(X)=h^0(X,K_X)=0.$$ Therefore, by the
Euler-Poincar\'e characteristic formula, we get
$$\chi(\sO_X)=1-h^{0,1}(X)+h^{0,2}(X)=0.$$ On the other hand, by the
Noether-Riemann-Roch formula,
$$\chi(\sO_X)=\frac{1}{12}(c_1^2(X)+c_2(X))=0,$$ we obtain
$$c_2(X)=-c_1^2(X).$$
Note also that $c_2(X)$ is the Euler characteristic $e(X)$ of $X$,
i.e.
$$c_2(X)=e(X)=2-2b_1(X)+b_2(X)=b_2(X)$$
and so $c_1^2(X)=-b_2(X)\leq 0$.  Suppose $X$ has a Levi-Civita
Ricci-flat Hermitian metric, then we have $c_1^{\mathrm{AC}}(X)=0$.
By Corollary \ref{top},  we have $c_1^2(X)=0$. Therefore $b_2(X)=0$.
It is well-known that, blowing-ups increase the second Betti number
at least by $1$, hence  we have $X=X_{\mathrm{min}}$. We complete
the proof of the statement that: if a compact complex surface admits
a Levi-Civita Ricci-flat metric, then it is a minimal surface.\\

There are three classes of surfaces of $\mathrm{VII}_0$:

 \bd

\item[$\bullet$] class $\mathrm{VII}_0$ surfaces with $b_2>0$;

\item[$\bullet$] Inoue surfaces: a class $\mathrm{VII}_0$ surface has $b_2=0$ and
contains no curves;

\item[$\bullet$] Hopf surfaces: its universal covering is $\C^2-\{0\}$, or
equivalently a class $\mathrm{VII}_0$ surface has $b_2=0$ and
contains a curve.

\ed

(1). A class $\mathrm{VII}_0$ surface $X$ with $b_2>0$ has no
Levi-Civita Ricci-flat metrics. Indeed, by a similar computation as
before, we know $c_1^2(X)=-b_2<0$ which contradicts to Corollary
\ref{top}.

(2). On an Inoue surface $X$, there is no Levi-Civita Ricci-flat
Hermitian metrics. This is essentially proved in
\cite[Theorem~7.2]{LY17}. For the reader's convenience, we include a
sketched proof here.
 It is well-known (\cite{In}) that an
Inoue surface is a quotient of $\H\times \C$ by a properly
discontinuous group of affine transformations where $\H$ is the
upper half-plane. There are three types of Inoue surfaces:

\bd \item[(A)] Inoue surfaces $S_M$. Let $M$ be a matrix in
$\mathrm{SL}_3(\Z)$ admitting one real eigenvalue $\alpha>1$ and two
complex conjugate eigenvalues $\beta\neq \bar\beta$. Let
$(a_1,a_2,a_3)$ be a real eigenvector of $M$ corresponding to
$\alpha$ and let $(b_1,b_2,b_3)$ be an eigenvector of $M$
corresponding to $\beta$. Then $X=S_M$ is the quotient of $\H\times
\C$ by the group of affine automorphisms generated by
$$g_0(w,z)=(\alpha w ,\beta z), \ \ \ \ g_i(w,z)=(w+a_i,z+b_i), \ \ i=1,2,3.$$

\item[(B)] Inoue surfaces $X=S^+_{N,p,q,r;t}$ are defined as the
quotient of $\H\times \C$ by the group of affine automorphisms
generated by
$$g_0(w, z)=(\alpha w, z+t),\ \ \ \ \ g_i(w,z)=(w+a_i, z+b_iw+c_i),\ \ i=1,2$$
$$g_3(w,z)=\left(w,z+\frac{b_1a_2-b_2a_1}{r}\right),$$
where $(a_1,a_2)$ and $(b_1,b_2)$ are the eigenvectors of some
matrix $N\in \mathrm{SL}_2(\Z)$ admitting real eigenvalues
$\alpha>1$, $\alpha^{-1}$. Moreover $t\in \C$ and $p, q, r (r\neq
0)$ are integers, and $(c_1,c_2)$ depends on $(a_i,b_i), p, q, r$.

\item[(C)] Inoue surfaces $X=S^-_{N,p,q,r;t}$  have unramified double cover which are Inoue surfaces of
type $S^+_{N,p,q,r;t}$.

\ed

\noindent
 Suppose--to the contrary--that  there exists a Levi-Civita Ricci-flat Hermitian
metric $\omega$ on the Inoue surface $X$.  Let
$\omega_f=e^{f}\omega$ be the Gauduchon metric
 in the conformal class of $\omega$, then by formula (\ref{Gauduchon}), the
total Chern scalar curvature of $\omega_f$ is \beq \int_X
s_{{f}}\cdot \omega_f^2=2\int_X\text{Ric}(\omega_f)\wedge
\omega_f=2\|\bp_f^*\omega_f\|_{\omega_f}^2\geq 0.\label{Inoue1}\eeq
We shall show that on each Inoue surface, there exists a smooth
Gauduchon metric with non-positive but not identically zero first
Chern-Ricci curvature. Indeed, it is easy to see that the metric
$h^{-1}=[\text{Im}(w)]^{-1}(dw\wedge dz)\ts (d\bar w\wedge d\bar z)$
(resp. $h^{-1}=[\text{Im}(w)]^{-2}(dw\wedge dz)\ts (d\bar w\wedge
d\bar z)$) is a globally defined Hermitian metric on the
anti-canonical bundle of $S_M$ (resp. $S^+_{N,p,q,r;t}$) (e.g.
\cite[Section~6]{DPS}). Hence, the Chern Ricci curvature of ${S_M}$
is
$$-\sq \p\bp\log h^{-1}=\sq\p\bp\log [\text{Im}(w)]=-\frac{\sq}{4}\frac{dw\wedge d\bar w}{[\text{Im}(w)]^2},$$
which also represents $c^{\mathrm{BC}}_1(X)$. By Theorem
\cite[Theorem~1.3]{STW}, there exists a Gauduchon metric $\omega_G$
with $$\text{Ric}(\omega_G)=-\frac{\sq}{4}\frac{dw\wedge d\bar
w}{[\text{Im}(w)]^2}\leq 0.$$ Hence, for any Gauduchon metric
$\omega$, one has
$$\int_X\text{Ric}(\omega)\wedge \omega=\int_X\text{Ric}(\omega_G)\wedge \omega<0$$
 which is a
contradiction to (\ref{Inoue1}). We can deduce similar
contradictions for $S^\pm_{N,p,q,r;t}$.

(3). A compact complex surface $X$ is called a Hopf surface if its
universal covering is analytically isomorphic to
$\C^2\setminus\{0\}$. It has been prove by Kodaira that its
fundamental group $\pi_1(X)$ is a finite extension of an infinite
cyclic group generated by a biholomorphic contraction which takes
the form \beq (z,w)\>(az,bw+\lambda z^m) \eeq where $a,b,\lambda \in
\C$, $|a|\geq |b|>1$, $m\in \N^*$ and $\lambda(a-b^m)=0$. Hence,
there are two different cases:

\bd \item[(I)] the Hopf surface $H_{a,b}$ of class $1$ if
$\lambda=0$;
\item[(II)] the Hopf surface $H_{a,b,\lambda,m}$ of class $0$ if $\lambda\neq
0$ and $a=b^m$. \ed

\noindent In the following, we consider the Hopf surface of class
$1$. Let $H_{a,b}=\C^{2}\setminus\{0\}/\sim$ where
$(z,w)\sim(az,bw)$ and $|a|\geq |b|> 1$. We set $k_1=\log|a|$ and
$k_2=\log|b|$. Define a real smooth function \beq
\Phi(z,w)=e^{\frac{k_1+k_2}{2\pi}\theta}\eeq where $\theta(z,w)$ is
a real smooth function defined by \beq
|z|^2e^{-\frac{k_1\theta}{\pi}}+|w|^2e^{-\frac{k_2\theta}{\pi}}=1.\label{0}\eeq
This is well-defined since for fixed $(z,w)$ the function
$t\>|z|^2|a|^t+|w|^2|b|^t$ is strictly increasing with image $\R_+$.
Let $\alpha=\frac{2k_1}{k_1+k_2}$ and so $1\leq\alpha<2$. Then the
key equation (\ref{0}) is equivalent to \beq
|z|^2\Phi^{-\alpha}+|w|^2\Phi^{\alpha-2}=1.\label{ke} \eeq It is
easy to see that
$$\theta(az,bw)=\theta(z,w)+2\pi, \qtq{and}\Phi(az,bw)=|a||b|\Phi(z,w).$$
We define a quantity \beq  \Delta=\alpha|z|^2
\Phi^{-\alpha}+(2-\alpha)|w|^2\Phi^{\alpha-2}.\label{delta} \eeq In
the next theorem, we construct precisely  Levi-Civita Ricci-flat
metrics on Hopf surfaces of class $1$. \eproof

\begin{theorem}\label{11}
On the Hopf surface $H_{a,b}$ of class $1$,  the Hermitian metric
\beq \omega=\Delta^3\left(\frac{1}{
\Phi}\sqrt{-1}\partial\bar\partial\Phi-\frac{1}{2}\sqrt{-1}\partial\bar\partial\log\Phi\right)
\label{LC} \eeq is  Levi-Civita Ricci-flat, i.e.
$\mathfrak{Ric}(\omega)=0$.
\end{theorem}

\bremark The proof of Theorem \ref{11} is carried out in the next
section. We should point out the construction follows from the ideas
in \cite[Theorem~6.2]{LY14} and \cite[Theorem~7.3]{LY17}. More
precisely, when $a=b$, we have $\alpha=1$, $\Delta=1$ and
$\Phi=|z|^2+|w|^2$. In this case, the Levi-Civita Ricc-flat metric
constructed in Theorem \ref{11} is exactly the same as the metrics
constructed in \cite[Theorem~6.2]{LY14} and
\cite[Theorem~7.3]{LY17}. \eremark

\vskip 2\baselineskip

\section{The construction of Levi-Civita Ricci-flat metrics on Hopf surfaces of type $1$}

In this section, we prove Theorem \ref{11}.

\blemma\label{aa} $\displaystyle |z|^2\Phi^{-\alpha}$ and
$|w|^2\Phi^{\alpha-2}$ are well-defined on $H_{a,b}$. \elemma
\bproof Indeed, $$ |az|^2\Phi^{-\alpha}(az,bw)=|a|^2
|a|^{-\alpha}|b|^{-\alpha} \cdot |z|^2\Phi^{-\alpha}(z,w) $$  and
$$|a|^2 |a|^{-\alpha}|b|^{-\alpha}=e^{k_1(2-\alpha)} e^{-k_2\alpha}=1.$$
Similarly, we can show $|w|^2\Phi^{2-\alpha}$ is well-defined on
$H_{a,b}$. \eproof

 \blemma\label{bb0} $\sq\p\bp\log\Phi$ has a
semi-positive matrix representation
\beq\frac{1}{\Delta^3\Phi^{2}}\left[\begin{array}{lcr}
(\alpha-2)^2|w|^2&\alpha(\alpha-2)\bar w z\\ \alpha(\alpha-2)\bar z
w& \alpha^2|z|^2
\end{array}\right],\label{222}\eeq
and $\sq\p\Phi\wedge \bp\Phi$ has a matrix representation \beq
\frac{1}{\Delta^2\Phi^{2\alpha-2}}\left[\begin{array}{lcr}
|z|^2&\bar w z\Phi^{2\alpha-2}\\ \bar z w \Phi^{2\alpha-2}&
|w|^2\Phi^{4\alpha-4}
\end{array}\right].\label{10}\eeq
 \elemma
\bproof See \cite[Appendix]{Y}. \eproof

 As motivated by \cite[Theorem~6.2]{LY14} and
\cite[Theorem~7.3]{LY17},  we consider the $(1,1)$-form \beq
\omega_\lambda=\frac{\sq
\p\bp\Phi}{\Phi}+\lambda\sq\p\bp\log\Phi.\label{lambda}\eeq It also
takes the form
$$\omega_\lambda=(1+\lambda)\sq\p\bp\log\Phi+\frac{ \sq \p\Phi\wedge \bp\Phi}{\Phi^2}$$
and it has the matrix representation \beq\left[\begin{array}{lcr}
\displaystyle{\frac{(1+\lambda)(\alpha-2)^2|w|^2}{\Delta^3\Phi^{2}}+\frac{|z|^2}{\Delta^2\Phi^{2\alpha}}}&\displaystyle{\frac{(1+\lambda)\alpha(\alpha-2)\bar w z}{\Delta^3\Phi^{2}}+\frac{\bar w z \Phi^{2\alpha-2}}{\Delta^2\Phi^{2\alpha}}}\\
\displaystyle{\frac{(1+\lambda)\alpha(\alpha-2)\bar z
w}{\Delta^3\Phi^{2}}+\frac{\bar z w
\Phi^{2\alpha-2}}{\Delta^2\Phi^{2\alpha}}}&
\displaystyle{\frac{(1+\lambda)\alpha^2|z|^2}{\Delta^3\Phi^{2}}+\frac{|w|^2\Phi^{4\alpha-4}}{\Delta^2\Phi^{2\alpha}}}
\end{array}\right].\eeq
Since $\det(\sq\p\bp\log\Phi)=\det(\sq \p\Phi\wedge \bp\Phi)=0$, the
determinant \be \det(\omega_\lambda)&=&
\frac{(1+\lambda)(\alpha-2)^2|w|^2}{\Delta^3\Phi^{2}}\cdot
\frac{|w|^2\Phi^{4\alpha-4}}{\Delta^2\Phi^{2\alpha}}+\frac{|z|^2}{\Delta^2\Phi^{2\alpha}}\cdot\frac{(1+\lambda)\alpha^2|z|^2}{\Delta^3\Phi^{2}}\\&&-\frac{(1+\lambda)\alpha(\alpha-2)\bar
w z}{\Delta^3\Phi^{2}}\cdot\frac{\bar z w
\Phi^{2\alpha-2}}{\Delta^2\Phi^{2\alpha}}-\frac{\bar w z
\Phi^{2\alpha-2}}{\Delta^2\Phi^{2\alpha}}\cdot
\frac{(1+\lambda)\alpha(\alpha-2)\bar z w}{\Delta^3\Phi^{2}}\\
&=&\frac{(1+\lambda)|w|^2(\alpha-2)}{\Delta^5}\left(\frac{(\alpha-2)|w|^2}{\Phi^{6-2\alpha}}-\frac{\alpha|z|^2}{\Phi^4}\right)\\&&+\frac{(1+\lambda)\alpha|z|^2}{\Delta^5}\left(\frac{\alpha|z|^2}{\Phi^{2\alpha+2}}-\frac{(\alpha-2)|w|^2}{\Phi^4}\right)\\
&=&\frac{1+\lambda}{\Delta^5}\cdot
\frac{1}{\Phi^{2+2\alpha}}\left(\alpha|z|^2+(2-\alpha)|w|^2\Phi^{2\alpha-2}\right)^2.
\ee By (\ref{delta}), we have
$$\alpha|z|^2+(2-\alpha)|w|^2\Phi^{2\alpha-2}=\Delta\Phi^\alpha$$
and so
 \beq
\det(\omega_\lambda)=\frac{1+\lambda}{\Delta^3\Phi^2}.
\label{det}\eeq It is easy to see from (\ref{lambda}) that, when
$\lambda>-1$, $$\omega_\lambda=\frac{\sq
\p\bp\Phi}{\Phi}+\lambda\sq\p\bp\log\Phi$$ is a Hermitian metric.
Let $\p^*$ and $\bp^*$ be the adjoint operators taken with respect
to the metric $\omega_\lambda$, and $\Lambda$ be the dual operator
of $\omega_\lambda\wedge\bullet$. \blemma We have \beq
\frac{\p\p^*\omega_\lambda+\bp\bp^*\omega_\lambda}{2}=\p\p^*\omega_\lambda=\bp\bp^*\omega_\lambda=\frac{\sq
\p\bp\log\Phi}{1+\lambda}.\label{tw}\eeq \elemma \bproof The metric
$\omega_\lambda$ has local matrix representation
$$ (h_{i\bar j})=
\begin{pmatrix}
\frac{1+\lambda}{\Phi}\Phi_{1\bar1}-\frac{\lambda\Phi_1\Phi_{\bar1}}{\Phi^2}&\frac{1+\lambda}{\Phi}\Phi_{1\bar2}-\frac{\lambda\Phi_1\Phi_{\bar2}}{\Phi^2}\\
\frac{1+\lambda}{\Phi}\Phi_{2\bar1}-\frac{\lambda\Phi_2\Phi_{\bar1}}{\Phi^2}&\frac{1+\lambda}{\Phi}\Phi_{2\bar2}-\frac{\lambda\Phi_2\Phi_{\bar2}}{\Phi^2}
\end{pmatrix}
$$ and its inverse matrix representation is $$
(h^{i\bar j})=\frac{\Phi^2\Delta^3}{1+\lambda}
\begin{pmatrix}
\frac{1+\lambda}{\Phi}\Phi_{2\bar2}-\frac{\lambda\Phi_2\Phi_{\bar2}}{\Phi^2}&-\frac{1+\lambda}{\Phi}\Phi_{1\bar2}+\frac{\lambda\Phi_1\Phi_{\bar2}}{\Phi^2}\\
-\frac{1+\lambda}{\Phi}\Phi_{2\bar1}+\frac{\lambda\Phi_2\Phi_{\bar1}}{\Phi^2}&\frac{1+\lambda}{\Phi}\Phi_{1\bar1}-\frac{\lambda\Phi_1\Phi_{\bar1}}{\Phi^2}
\end{pmatrix}.
$$
By Lemma \ref{adjoint}, we have
$$\p^*\omega_\lambda=-\sq\Lambda\bp\omega_\lambda=-2\sqrt{-1}(\Gamma_{\lambda})^k_{\bar ik}d\bar z^i.$$
A straightforward computation shows that \be \frac{\p
h_{j\bar\ell}}{\p \bar z^i}&=&\frac{\p}{\p\bar
z^i}\left(\frac{1+\lambda}{\Phi}\Phi_{j\bar\ell}-\frac{\lambda}{\Phi^2}\Phi_j\Phi_{\bar\ell}\right)\\
&=&-\frac{1+\lambda}{\Phi^2}\Phi_{\bar
i}\Phi_{j\bar\ell}+\frac{1+\lambda}{\Phi}\Phi_{j\bar\ell\bar
i}+\frac{2\lambda}{\Phi^3}\Phi_{\bar
i}\Phi_j\Phi_{\bar\ell}-\frac{\lambda}{\Phi^2}\Phi_{j\bar
i}\Phi_{\bar\ell}-\frac{\lambda}{\Phi^2}\Phi_{j}\Phi_{\bar\ell\bar
i}.\ee Hence, we have
\begin{align*}
(\Gamma_{\lambda})^k_{\bar ij}=&\frac{1}{2}h^{k\bar \ell}\left(\frac{\partial h_{j\bar  \ell}}{\partial\bar z^i}-\frac{\partial h_{j\bar i}}{\partial\bar z^\ell}\right)\\
=&\frac{1}{2}h^{k\bar \ell}\left( -\frac{1}{\Phi^2}\Phi_{j\bar
\ell}\Phi_{\bar i} +\frac{1}{\Phi^2}\Phi_{j\bar i}\Phi_{\bar \ell}
\right)
\end{align*}
and  \be
\partial^*\omega_{\lambda}&=&-2\sqrt{-1}\left(\Gamma_{\lambda}\right)^k_{\bar ik}d\bar z^i=-\sqrt{-1}h^{k\bar \ell}\frac{1}{\Phi^2}\left(-\Phi_{k\bar \ell}\Phi_{\bar i}
+\Phi_{k\bar i}\Phi_{\bar \ell}\right)d\bar z^i\\
&=& -\sqrt{-1}h^{k\bar \ell}\frac{1}{\Phi^2}\left(-\Phi_{\bar
i}\left(\frac{\Phi h_{k\bar \ell}}{1+\lambda}+\frac{\lambda
\Phi_{k}\Phi_{\bar \ell}}{(1+\lambda)\Phi}\right)
+\Phi_{\bar \ell}\left(\frac{\Phi h_{k\bar i}}{1+\lambda}+\frac{\lambda \Phi_{k}\Phi_{\bar i}}{(1+\lambda)\Phi}\right)\right)\\
&=&-\sqrt{-1}h^{k\bar \ell}\frac{1}{\Phi^2}\left(\frac{-\Phi_{\bar
i}\Phi h_{k\bar \ell}}{1+\lambda}
+\frac{\Phi_{\bar \ell}\Phi h_{k\bar i}}{1+\lambda}\right)\\
&=& \frac{\Phi_{\bar i}}{\Phi}\frac{\sqrt{-1}}{1+\lambda}d\bar
z^i=\frac{\sq}{1+\lambda} \bp\log\Phi. \ee Therefore, we get
(\ref{tw}). \eproof

\noindent By formulas (\ref{key}), (\ref{det}) and (\ref{tw}), we obtain \be  \mathfrak{Ric}(\omega_\lambda)&=&-\sqrt{-1}\partial\bar\partial\det (\omega_{\lambda})-\frac{\partial\partial^*\omega_{\lambda}+\bar\partial\bar\partial^*\omega_{\lambda}}{2}\\
&=&\left(2-\frac{1}{1+\lambda}\right)\sqrt{-1}\partial\bar\partial\log\Phi+3\sqrt{-1}\partial\bar\partial\log\Delta.
 \ee
In particular, we can take $\lambda=-\frac{1}{2}$, and obtain \beq
\mathfrak{Ric}(\omega_{-\frac{1}{2}})=3\sqrt{-1}\partial\bar\partial\log\Delta.
\label{half}\eeq

 \btheorem Let
$$\omega=\Delta^3\omega_{-\frac{1}{2}}=\Delta^3\left(\frac{\sqrt{-1}}{\Phi}\partial\bar\partial\Phi  -\frac{1}{2}\sqrt{-1}\partial\bar\partial\log\Phi\right).$$
Then we have \beq \mathfrak{Ric}(\omega)=0. \eeq\etheorem \bproof By
Lemma \ref{conformalchange} and formula (\ref{half}), we have \beq
\mathfrak{Ric}(\omega)=\mathfrak{Ric}(\omega_{-\frac{1}{2}})-\sq
\p\bp\log \Delta^3=0. \eeq

\eproof

%


\begin{thebibliography}{99}

\bibitem{AT13} Angella, D.; Tomassini, A. On the $\p\bp$-lemma and Bott-Chern cohomology. Invent. Math. \textbf{192} (2013), no. \textbf{1}, 71--81.


\bibitem{BHPV} Barth, W.; Hulek, K.;  Peters, C.;  Van de Ven, A. \textit{Compact complex surfaces}.  Ergebnisse der Mathematik und ihrer Grenzgebiete. 3. Folge. A Series of Modern Surveys in Mathematics. Springer-Verlag, Berlin, 2004.

\bibitem{CRS} Chiose, I.; Rasdeaconu, R.; Suvaina, I. Balanced metrics on uniruled
manifolds.
 \href{http://arxiv.org/abs/1408.4769}{arXiv:1408.4769}. To appear
 in Comm. Anal. Geom.




\bibitem{DPS}  Demailly, J.-P.; Peternell, T.; Schneider, M. Compact complex manifolds with numerically effective tangent bundles. J. Algebraic Geom. \textbf{3} (1994), no. 2, 295--345.


\bibitem{FTWZ} Fang, S.-W.;  Tosatti, V.;  Weinkove, B. and  Zheng, T.  Inoue surfaces and the Chern-Ricci
flow.  J. Funct. Anal. \textbf{271} (2016),  3162--3185.

 \bibitem{F3} Fu, J.-X.  On non-K\"ahler Calabi-Yau threefolds with balanced metrics. \emph{Proceedings of the International Congress of Mathematicians.} Volume II, 705--716, Hindustan Book Agency, New Delhi, 2010.
\bibitem{FY} Fu, J.- X.; Yau, S.-T. The theory of superstring with flux on non-K\"ahler manifolds and the complex Monge-Amp\'ere equation. J. Differential Geom. \textbf{78} (2008), no. 3, 369--428.
\bibitem{FLY}Fu, J.- X; Li, J.; Yau, S.-T. Constructing balanced metrics on some families of non-K\"ahler Calabi-Yau threefolds.  J. Differential Geom. \textbf{90} (2012), no. 1,81--129.

\bibitem{Ga2}  Gauduchon, P. { Fibr\'es hermitiens \`a endomorphisme de Ricci non-n\'egatif}, Bull. Soc. Math. France {\textbf {105}}  1977), 113--140.






\bibitem{In} Inoue, M. On surfaces of type $\mathrm{VII_0}$, Inv. Math., \textbf{24} (1974), 269--310.


\bibitem{L} Li, Yi. A priori estimates for Donaldson’s equation over compact Hermitian manifolds. Calc. Var. Partial Differential Equations, \textbf{50}(2014), no. 3-4, 867--882.

\bibitem{LY12} Liu, K.-F.; Yang, X.-K. Geometry of Hermitian manifolds.  Internat. J. Math. \textbf{23} (2012) 40pp.

\bibitem{LY14} Liu, K.-F.; Yang, X.-K. Ricci curvatures on Hermitian manifolds.{Trans. Amer. Math. Soc.}  \textbf{369} (2017),
5157--5196.

\bibitem{LY17} Liu, K.-F.; Yang, X.-K. , Minimal complex surface with
Levi-Civita Ricci flat metrics.\\
\href{https://arxiv.org/abs/1706.01219}{arXiv:1706.01219} \emph{To
appear in Acta. Math. Sinica.} In Memory of Professor Lu Qi-Keng.


\bibitem{ST2} Streets, J.; Tian, G. {A parabolic flow of pluriclosed metrics}, Int. Math. Res. Not.  2010, no. 16, 3101--3133.
\bibitem{ST3} Streets, J.; Tian, G. {Regularity results for pluriclosed flow}, Geom. Topol. \textbf{17} (2013), no. 4, 2389--2429.


\bibitem{STW} Sz\'{e}kelyhidi, G.;  Tosatti, V.; Weinkove, B. Gauduchon metrics with prescribed volume
form.  Acta Math. \textbf{219} (2017), no.1, 181--211.

\bibitem{Tel} Teleman, A. The pseudo-effective cone of a non-K\"ahlerian surface
and applications. Math. Ann. \textbf{335}(2006), 965--989.

\bibitem{T} Tosatti, V. Non-K\"ahler Calabi-Yau manifolds.  Contemp. Math. \textbf{644} (2015), 261--277.

\bibitem{TW} Tosatti, V.; Weinkove, B. On the evolution of a Hermitian metric by its Chern-Ricci
form. J. Differential Geom. \textbf{99} (2015), no.1, 125--163.

\bibitem{TW2} Tosatti, V.; Weinkove, B. The
Monge-Amp\`{e}re equation for $(n-1)$-plurisubharmonic functions on
a compact K\"ahler manifold. J. Amer. Math. Soc. \text{30} (2017),
no.2, 311-346.



\bibitem{Yang16}Yang, X.-K. Hermitian manifolds with semi-positive
holomorphic sectional curvature.  {Math. Res. Lett.} \textbf{23}
(2016), no.3, 939--952.

\bibitem{Yang17A} Yang, X.-K.  The Chern-Ricci
flow and holomorphic bisectional curvature.
 Sci. China Math. \textbf{59} (2016), 2199-2204.
\bibitem{Y} Yang, X.-K.  Big vector bundles and compact complex manifolds with semi-positive tangent
bundle. Math. Ann. \textbf{267}(2017), no.1, 251--282.
\bibitem{Yang17} Yang, X.-K. Scalar curvature on compact complex
manifolds.
\href{https://arxiv.org/abs/1705.02672}{arXiv:1705.02672}. To appear
in Trans. AMS.

\bibitem{Yau2} Yau, S.-T. On the curvature of compact Hermitian manifolds. Invent. Math. \textbf{25} (1974), 213--239.
\bibitem{Yau78} Yau, S.-T.  {On the Ricci curvature of a compact K\"ahler manifold and the complex Monge-Amp\`ere equation, I},
Comm. Pure Appl. Math. {\textbf{ 31}} (1978),  339--411.

\end{thebibliography}
\end{document}